\documentclass[amscd]{amsart}

\newtheorem{Thm}{Theorem}
\newtheorem{Cor}{Corollary}
\newtheorem{Lem}{Lemma}
\newtheorem{Prop}{Proposition}

\newtheorem{Claim}{Claim}
\newtheorem{Conj}{Conjecture}

\theoremstyle{remark}
\newtheorem{Rem}{Remark}

\newsymbol\ltimes 226E
\newsymbol\rtimes 226F
\newsymbol\boxtimes 1202
\newsymbol\twoheadrightarrow 1310
\newsymbol\rightarrowtail 131A

\def\square{\hbox{\vrule\vbox{\hrule\phantom{o}\hrule}\vrule}}

\newcommand{\cal}{\mathcal}

\newcommand{\Fl}{{{\cal F}\ell}}

\newcommand{\isol}{{\widetilde \longleftarrow}}
\newcommand{\isor}{{\widetilde \longrightarrow}}

\newcommand{\<}{\langle}
\renewcommand{\>}{\rangle}

\def\square{\hbox{\vrule\vbox{\hrule\phantom{o}\hrule}\vrule}}
\newcommand{\epf}{\square}

\newcommand{\CC}{{\cal C}}

\newcommand{\EE}{{\cal E}}

\newcommand{\fD}{{\mathfrak D}}

\newcommand{\bw}{{\bf w}}

\newcommand{\wti}{\tilde w}
\newcommand{\N}{{\cal N}}

\newcommand{\Nt}{{\tilde{\cal N}}}

\newcommand{\fS}{{\frak S}}

\newcommand{\Oe}{{\frak O}}

\newcommand{\g}{{\frak g}}
\newcommand{\fb}{{\frak b}}
\newcommand{\fp}{{\frak p}}
\newcommand{\gh}{\widehat{\frak g}}

\newcommand{\Lg}{{ \frak g\check{\ }}}

\newcommand{\LG}{{ G\check{\ }}}

\renewcommand{\O}{{\cal O}}

\newcommand{\F}{{\cal F}}
\newcommand{\OO}{{\cal O}}
\newcommand{\GG}{{\cal G}}

\newcommand{\BB}{{\cal B}}

\newcommand{\Zet}{{\Bbb Z}}
\newcommand{\Ce}{{\Bbb C}}

\newcommand{\GO}{{\bf G_O}}

\newcommand{\Pone}{{\mathbb P}^1}

\newcommand{\Ga}{\Gamma}

\newcommand{\new}{{\mathrm{new}}}

\newcommand{\Opn}{Op_{\mathrm{nilp}}}

\newcommand{\Opr}{Op_{\mathrm{rs}}}

\newcommand{\Opi}{Op_{\mathrm{int}}}

\newcommand{\Hom}{{\mathrm{Hom}}}

\newcommand{\la}{\lambda}

\newcommand{\La}{\Lambda}

\newcommand{\Ltimes}{\overset{\rm L}{\times}}

\newcommand{\gt}{\tilde{\g}}

\newcommand{\tii}{\widetilde}

\newcommand{\balp}{\aleph}

\newcommand{\st}{\star}

\newcommand{\tmod}{\mbox{--}{\mathrm{mod}}}

\title[Critical level and noncommutative Springer]
{Highest weight modules at the critical level and noncommutative Springer resolution\footnote{AMS
 S\lowercase{ubject} C\lowercase{lassification} (2010):
17B67, 22E57, 14D24,
20G42.}\footnote{T\lowercase{he first named author is supported by} NSF
\lowercase{grants} DMS-0854764  \lowercase{and} DMS-1102434.}
}

\author{
Roman Bezrukavnikov
}
\address{\small
Department of Mathematics, Massachusetts Institute of Technology, 77 Massachusetts ave.,
Cambridge, MA 02139, USA
}
\email{
bezrukav@math.mit.edu
}

\author{
Qian Lin
}
\address{\small
Oracle Hardware and Software Engineering,
400 Oracle Parkway, Redwood City, CA 94065
}
\email{
qianlin88@gmail.com
}
\begin{document}

\maketitle

\begin{abstract}
In \cite{BM} a certain non-commutative algebra $A$ was defined starting from a semi-simple algebraic group, so that the derived category of $A$-modules is equivalent to the derived category of coherent sheaves on the Springer (or Grothendieck-Springer) resolution.

Let $\Lg$ be the Langlands dual Lie algebra and let
$\gh$ be the corresponding affine Lie algebra,
i.e. $\gh$ is a central extension of $\Lg\otimes \Ce((t))$.

Using results of Frenkel and Gaitsgory
we show that the category of $\gh$ modules at the critical level which are Iwahori integrable
and have a fixed central character,  is equivalent to the category of modules over a
central reduction of $A$. This implies that numerics of Iwahori integrable modules at the critical level is governed by the canonical basis in the $K$-group of a Springer fiber, which was conjecturally described by Lusztig \cite{Kth2} and constructed in \cite{BM}.

\end{abstract}

\tableofcontents

\section{Introduction}
\subsection{Modules at the critical level}
Category $O$ of highest weight modules over a semi-simple Lie algebra with a fixed central character is a classical object of study in representation theory; Kazhdan-Lusztig conjectures (proved by Beilinson--Bernstein and Brylinksi--Kashiwara) assert that numerics of
such modules is governed by the canonical  basis in the Hecke algebra.
The subject of this paper is an analogue of that result for
modules over an affine Lie algebra at the critical level. We show that the category of such modules
is governed by the canonical bases in the Grothendieck group (or homology) of Springer fibers. This basis was described conjecturally by Lusztig \cite{Kth2} and its existence was established in \cite{BM}.
 The original motivation
for \cite{BM} came from representation theory of Lie algebras in positive characteristic; it turns out that the same generalization of Kazhdan-Lusztig theory controls highest weight modules at the critical level.

Let $G$ be a semi-simple algebraic group 
 over $\Ce$ with Lie algebra $\g$.
Let $\gh$ be the affine Lie algebra corresponding to the {\em Langlands dual}
Lie algebra $\Lg$. Thus $\gh$ is a central extension of the loop algebra
$$0\to \Ce^r\to \gh\to \Lg\otimes \Ce((t))\to 0,$$
where $r$ is the number of simple summands in $\g$.

Let $U_{crit}\gh$ denote the quotient of the enveloping algebra $U\gh$ at the critical value of the central charge. 

Let $U_{crit}\gh\tmod^{I^0}$ denote the category of Iwahori monodromic $U_{crit}\gh$ modules in the sense of \cite{FG1}. Recall that by the result of Feigin and Frenkel
\cite{FF} a continuous  $U_{crit}\gh$ module
carries a canonical commuting action of the topological ring $\OO(Op)$ of functions on the space $Op$ of $G$-opers on the formal punctured disc. In particular, an irreducible module has a central character which corresponds to such an oper. For an irreducible module $L\in U_{crit}\gh\tmod^{I^0}$
the oper necessarily has a regular singularity and a nilpotent residue.

Fix a nilpotent element $e\in \g$ and a nilpotent oper $\Oe$ with residue $e$ (thus,
the underlying connection is isomorphic to $\nabla=d+e\frac{dt}{t}$ where $t$ is a coordinate on the formal disc).
We let $U_{crit}\gh\tmod^{I^0}_{\Oe}$ be the full subcategory in $U_{crit}\gh\tmod^{I^0}$
consisting of finite length modules where $\OO(Op)$ acts through the character corresponding to $\Oe$.

\subsection{Noncommutative Springer resolution}
We now introduce another abelian category associated to the nilpotent element $e$.
Let
$\BB=G/B$ be the flag variety of $G$
thought of as the variety of Borel subalgebras in $\g$; let  $\gt=\{(x,\fb)\ |\fb \in \BB, x\in \fb \}\overset{\pi}{\longrightarrow} \g$ be the Grothendieck-Springer map $\pi:(\fb,x)\mapsto x$.


In \cite{BM} a certain non-commutative algebra $A$, well defined up to a Morita equivalence was introduced. The algebra comes equipped with an equivalence of triangulated categories $D^b(A\tmod^{fg})\cong D^b(Coh(\gt))$; by
$A\tmod^{fg}$ we denote the category of finitely generated $A$-modules.

The center of $A$ is identified with the algebra $\O(\g)$ of regular functions on $\g$. For $e\in \g$ we let $A_e$ denote the corresponding central reduction of $A$.

The results of \cite{BM} provide a canonical isomorphism of Grothendieck groups
$K^0(A_e\tmod^{fg})\cong K^0(Coh(\pi^{-1}(e)))$ sending the classes of irreducible modules to elements of the {\em canonical basis}, i.e. the
unique (up to signs) basis satisfying the axioms of \cite{Kth2}.

The next statement conjectured in \cite[Conjecture 1.7.2]{BM} is the main result of the present note.

\begin{Thm}\label{main}
There exists a canonical equivalence of abelian categories
$$A_e\tmod^{fg}\cong U_{crit}\gh\tmod^{I^0}_{\Oe}.$$
\end{Thm}

In fact, the equivalence of derived categories follows by comparing the result of
 \cite{FG} which identifies $D^b(U_{crit}\gh\tmod^{I^0}_{\Oe})$ with $DGCoh(\BB_e)$ with
 that of
 \cite{BM}, \cite{Ri} which identify $D^b(A_e\tmod^{fg})$ with the same category
 of coherent sheaves.
Here  $DGCoh(\BB_e)$  denotes the derived category of coherent sheaves on the DG-scheme
$\BB_e=\{e\} \Ltimes_{\g} \gt$. (Notice that $K^0(DGCoh(\BB_e))=K^0(Coh(\pi^{-1}(e)))$
since the Grothendieck group of coherent complexes on a DG-scheme is identified with the Grothendieck group of coherent sheaves on the underlying scheme.)

Our job in the present note is to show that the resulting equivalence
$$D^b(A_e\tmod^{fg})\cong D^b(U_{crit}\gh\tmod^{I^0}_{\Oe})$$ induces an equivalence of abelian categories, i.e. that it is $t$-exact with respect to the natural $t$-structures.
This will be done using characterizations of the $t$-structure on $DGCoh(\BB_e)$
coming from the two equivalences with derived categories of modules, appearing,  respectively, in \cite{BM}, \cite{FG}.

According to \cite[Conjecture 1.7.1]{BM}, the category of $A_e\tmod^{fg}$ is equivalent to a category of modules over the Kac - De Concini quantum group at a root of unity. Thus, together with the present result, that Conjecture implies an equivalence between modules over the affine Lie algebra and quantum group modules at a root of unity. Another equivalence of this sort has been established in the celebrated work by Kazhdan and Lusztig \cite{KL}.

We also expect that  when the nilpotent $e$ is of principal
Levi type (a generalization of) our result can be used to derive character formulas for irreducible highest weight modules in terms of parabolic periodic Kazhdan-Lusztig polynomials; we plan to develop this application in a future work.

\bigskip

The rest of the text is structured as follows. In section \ref{NCS}
we recall the needed properties of the noncommutative Springer resolution including
 a characterization of the corresponding $t$-structure on the derived categories of coherent sheaves.
Section \ref{affl} is devoted to constructible sheaves on affine flag variety of the dual group. We state a description of the subcategory of complexes
equivariant with respect to the radical of the Iwahori subgroup $I^0$
 in terms of coherent sheaves
on Steinberg variety of $G$, to appear in \cite{tobe}.
A technical result
about the $t$-structure on the category of Iwahori-Whittaker sheaves
appearing in Proposition \ref{iffIW}
is the key statement providing a link between the description of the $t$-structure
by Frenkel-Gaitsgory \cite{FG} to our formalism of braid positive $t$-structures.
In section \ref{CRIT} we quote the result of \cite{FG}
and argue, in subsection \ref{COMP}, that the two characterizations are compatible, which yields Theorem \ref{main}.


\subsection{Conventions and notations}
Let $W_{aff}$ denote the semi-direct product of the Weyl group $W$ by the weight lattice $\Lambda$ of $G$.
Thus $W_{aff}$ is an extended affine Weyl group
corresponding (in the Bourbaki terminology) to  the {\em dual group} $\LG$.
 Let
$\ell$ denote the length function on $W_{aff}$.
 Let $B_{aff}$ be the  corresponding extended
  affine braid group
and $B_{aff}^+\subset B_{aff}$ be the semigroup of positive braids, i.e. the semigroup consisting
of products of the Coxeter generators (but not their inverses) and length zero elements.
Thus $B_{aff}$ surjects onto  $W_{aff}$.
  We have a section of the map $B_{aff}\to W_{aff}$ sending an element $w\in W_{aff}$ to its minimal length preimage $\wti\in B_{aff}^+$. The elements
$\wti$ generate $B_{aff}$ subject to the relation $\tii{w_1w_2}=\wti_1\wti_2$ provided that
$\ell(w_1w_2)=\ell(w_1)+\ell(w_2)$.
Let $\Sigma$ (respectively, $\Sigma_{aff}$) be the set of vertices of  Dynkin diagram (respectively,
affine Dynkin diagram) of $\Lg$. For dominant weights $\la,\mu \in \Lambda^+\subset \Lambda$
we have $\tii{\la} \tii{\mu}=\tii{\la+\mu}$, thus we have a homomorphism $\Lambda\to B_{aff}$
sending $\la\in \Lambda^+$ to $\tii{\la}$, we denote this homomorphism by $\la\mapsto \theta_\la$.

For a set $S$ of objects in a triangulated category $\CC$ we let $\<  S\>$ denote
the full subcategory in $\CC$ generated by $S$ under extensions and direct summands.

\subsection{Acknowledgements} We thank Dennis Gaitsgory for helpful discussions.

\section{Noncommutative Springer resolution}\label{NCS}
In this section we summarize the results of \cite{BM}, \cite{BR} (see also \cite{ICM}).

We will use the language of DG-schemes, see \cite{BR} for the summary of necessary elementary facts.
We only use DG-schemes explicitly presented as fiber products of ordinary schemes, so we do not require the subtler aspects of the theory discussed in the current literature on the subject.

 The concept of a {\em geometric action} of a group
on a scheme $X$ over a scheme $Y$ (where a map $X\to Y$ is fixed) is introduced in \cite{BM}, \cite{BR}.
We do not recall the definition in detail, but we mention that a geometric action induces
a usual action on the derived category $DGCoh(X\Ltimes_Y S)$ for any scheme $S$ mapping to $Y$. Here
$X\Ltimes_Y S$ is the derived fiber product and $DGCoh$ denotes the triangulated category of
sheaves of coherent $DG$-modules over the structure sheaf. In the case when higher Tor sheaves
$Tor_i^{\O(Y)}(\O(X), \O(S))$, $i>0$, vanish, this reduces to the usual fiber product
$X\times _Y S$ and we have $DGCoh(X\Ltimes_Y S)=D^b(Coh(X\times _Y S))$.
For varying $S$, the actions are compatible with
pull-back and push-forward functors.

Recall that $\pi:\gt \to \g$ is the Grothendieck-Springer map.
In \cite{BR} a geometric action of $B_{aff}$ on $\gt$ over $\g$ is constructed.

For a quasi-projective scheme $S$ of finite type over $\Ce$ with
a fixed map to $\g$
set $\gt_S=\gt\Ltimes_\g S$.

We let $\balp$ denote the geometric action and $\balp_S$ the corresponding action of $B_{aff}$ on $DGCoh(\gt_S)$. The action $\balp_S$ can be described as follows.

For $\la\in \La$ let $\O_\BB(\la)$ be the corresponding line bundle on the flag variety, and $\OO_{\gt_S}$
be its pull-back to $\gt_S$.

For $\alpha\in \Sigma$ let ${\mathcal P}_\alpha$ be the corresponding partial flag variety
thought of as the variety of parabolic subalgebras in $\g$ belonging to a fixed conjugacy class. Let
 $\gt_\alpha=\{ (x,\fp)\ |\ \fp\in {\mathcal P}_\alpha, x\in \fp\}$,  and let $\Ga^\alpha$ denote the component of $\gt\times _{\gt_\alpha}\gt$ different from the diagonal. Let $\Ga^\alpha_S=\Ga^\alpha\Ltimes_\g S$. Let $pr_1,\, pr_2:\Ga^\alpha_S\to \gt_S$ be the projections.

Then we have:
$$\aleph_S(\tii s_\alpha): \F\to pr_{2*}^\alpha pr_{1}^{\alpha*}(\F),\ \ \ \ \alpha\in \Sigma$$
$$\aleph_S(\theta_\lambda): \F \to \F\otimes_{\O_{\gt_S}}\O_{\gt_S}(\lambda).$$

We say\footnote{This terminology differs slightly from that of \cite{BM} --  a normalized braid positive $t$-structure was called an exotic $t$-structure in {\em loc. cit.}}
that a $t$-structure $\tau$ on $DGCoh(\gt_S)$ is
 {\em braid positive} if $\aleph_S(\tii{s_\alpha})$, $\alpha\in \Sigma_{aff}$
  is right exact
 i.e. it sends $DGCoh(\gt_S)^{\leq 0}_\tau$ to itself.  Notice that this
definition involves the action of $\tii{s_\alpha}$ for all $\alpha \in \Sigma_{aff}$; in particular,
for $\alpha\not \in \Sigma$ this action is not given by an explicit correspondence
(though it can be expressed as a composition of correspondences used defining
the action of $\tii{s_\alpha}^{\pm 1}$,
$\alpha\in \Sigma$ and $\theta_\lambda$).

We will say that such a $t$-structure
is {\em normalized} if the direct image functor $R\pi_{S*}: DGCoh(\gt_S)\to D^b(Coh(S))$
is $t$-exact where the target category is equipped with the tautological $t$-structure.

The following was established in \cite{BM}.

\begin{Thm}\label{brpex}
a) For any $S$ a normalized braid positive $t$-structure exists and it is unique.
It satisfies:

$\F\in D^{\leq 0} $ iff  $pr_{S*} (b(F))\in D^{\leq 0}(Coh(S))$ for all $b\in B_{aff}^+$,

 $\F\in D^{\geq 0} $ iff $pr_{S*} (b^{-1}(F))\in D^{\geq 0}(Coh(S))$ for all $b\in B_{aff}^+$.

b) There exists a finite locally free $\O(\g)$ algebra $A$ such that for any $S$ as above
there is an equivalence $D^b(A_S\tmod ^{fg})\cong DGCoh(\gt_S)$ sending the tautological $t$-structure
on the LHS to the normalized braid positive $t$-structure on the RHS. Here $A_S=A\otimes _{\O(\g)}\O(S)$ and $A_S\tmod ^{fg}$ is the category of finitely generated $A_S$-modules.

\end{Thm}

\begin{Rem}
It is easy to see that the properties stated in part (b) of the Theorem characterize the algebra $A$ appearing there uniquely up to a Morita equivalence. (In fact, the part of the statement pertaining to the absolute case $S=\g$ is sufficient to characterize $A$). For notational convenience we fix a representative $A$ of the Morita equivalence class.
\end{Rem}

\begin{Rem}
The Theorem was stated in \cite{BM} under the additional assumption of Tor vanishing,
when $\gt_S$ can be considered as an ordinary scheme rather than a DG-scheme, and only for affine $S$. However, the proof carries over to the case of arbitrary base change involving DG-schemes, given the foundational material in
\cite[\S 1]{Ri}, \cite{BR}.
\end{Rem}

\begin{Rem}
The characterization of a normalized braid positive $t$-structure involves only the action
of Coxeter generators $\tii s_\alpha$ which generate the semi-group $B_{aff}^+$
if $G$ is adjoint but not in general. However, elements of $B_{aff}^+$ act by right exact functors
for any $G$, see \cite[Remark 1.5.2]{BM}. In particular, the subgroup $\Omega$ of length zero elements
in $W_{aff}$ acts by $t$-exact automorphisms, i.e. it acts by automorphisms of the corresponding abelian
heart. Notice that $\Omega$ acts on $\gh$ by outer automorphisms coming from automorphisms of the affine
Dynkin diagram.  Thus if an oper $\Oe$
is $\Omega$-invariant we get an action of $\Omega$ on $U_{crit}\gh\tmod^{I^0}_{\Oe}$.
It is natural to conjecture
that for such an oper the equivalence of Theorem \ref{main} is compatible with the action
of $\Omega$.
\end{Rem}


\subsection{Base change to a point and canonical bases} We now turn to the particular case
when $S=\{e\}$ is a point. We assume for simplicity that $e\in \N$.

Then we get a finite dimensional algebra $A_e$ together with the equivalence
\begin{equation}\label{ee}
D^b(A_e\tmod^{fg})\cong DGCoh(\gt_e).
\end{equation}

The reduced variety of the DG-scheme $\gt_e$ is the {\em Springer fiber} $\BB_e=\pi^{-1}(e)$.
It follows that the Grothendieck group $K^0(DGCoh(\gt_e))$ is isomorphic to $K^0(Coh(\BB_e))=:K_0(\BB_e)$.
The equivalence \eqref{ee} induces an isomorphism
$$K^0(A_e\tmod^{fg})\cong K^0(DGCoh(\gt_e)) = K^0(Coh(\BB_e)).$$

Since $A_e$ is a finite dimensional algebra over $k$, the group $K^0(A_e\tmod^{fg})$ is a free abelian group with a basis formed by the classes of irreducible $A_e$-modules.

The following is a restatement of the main result of \cite{BM}.

\begin{Thm}
The basis in $K_0(\BB_e)$ corresponding to the basis of irreducible $A_e$-modules under the above isomorphisms is the canonical basis,
i.e. it is characterized (uniquely up to signs) by the axioms of \cite{Kth2}.
\end{Thm}

\subsection{The equivariant version} Let $H$ be a reductive group with a homomorphism $H\to G$,
and assume that $H$ acts on $S$ so that the map $S\to \g$ is $H$-equivariant.
Then we get (see \cite[1.6.6]{BM})
\begin{equation}\label{equiveqnt}
DGCoh^H(\gt_S)\cong D^b(A_S\tmod_{coh}^H),
\end{equation}
where $A_S\tmod_{coh}^H$ denotes the category $H$-equivariant finitely generated $A_S$-modules.

Below we will apply it in the case $S=\Nt$, $H=G$.

\section{Perverse sheaves on the affine flag variety}\label{affl}

\subsection{Affine flag variety and categories of constructible sheaves}
Along with the Lie algebra $\Lg((t))$ we consider the group ind-scheme $\LG((t))$ and its group subschemes
$I^0\subset I\subset \GO$, where $I$ is the Iwahori subgroup and $I^0$ is its pro-unipotent radical
and $\GO=\LG[[t]]$ is the subgroup of regular loops into $\LG$. Let also $I_-\subset \GO$ be an opposite Iwahori and $I^0_-$ be its pro-unipotent radical.
Let $\Fl$ be the affine flag variety for the group $\LG$; thus $\Fl=\LG((t))/I$ is an ind-projective
ind-scheme.

We will consider the following full subcategories in the derived category $D(\Fl)$
of constructible sheaves
on $\Fl$: the category $D_{I^0}(\Fl)$ of complexes equivariant with respect to $I^0$ and $D_{IW}$ of complexes equivariant with respect to a non-degenerate character of $I_-^0$, see \cite{AB} for details.
The functors of forgetting the equivariance $D_{I^0}(\Fl)\to D(\Fl)$, $D_{IW}\to D(\Fl)$ are full embeddings since the group schemes $I^0$, $I_-^0$ are pro-unipotent.

Let $Perv_{I^0}(\Fl)\subset D_{I^0}(\Fl)$, $Perv_{IW}(\Fl)\subset D_{IW}(\Fl)$ be the full subcategories
of perverse sheaves.
It is known that there are natural equivalence
$D^b(Perv_{I^0}(\Fl))\cong D_{I^0}(\Fl)$, $D^b(Perv_{IW}(\Fl))\cong D_{IW}(\Fl)$ (see e.g. \cite{AB}).

We will also need the Iwahori equivariant derived category $D_I(\Fl)$ (which, in contrast with
the categories $D_{I^0}(\Fl)$, $D_{IW}$ is not equivalent to the derived category of the abelian subcategory of perverse sheaves $Perv_I(\Fl)$). The category $D_I(\Fl)$ carries a monoidal structure provided by convolution which will be denoted by $\st$. This monoidal category acts on $D(\Fl)$ by convolution on the right, which will also be denoted by $\st$; the action preserves the subcategories
$D_{I^0}$, $D_{IW}$.

The orbits of $I$ on $\Fl$ are indexed by $W_{aff}$; for $w\in W_{aff}$ let $\Fl_w$ denote the corresponding orbit and $j_w:\Fl_w\to \Fl$ be the embedding. We abbreviate $j_{w*}=j_{w*}(\underline{\Ce}[\dim \Fl_w])$, $j_{w!}=j_{w!}(\underline{\Ce}[\dim \Fl_w])$; thus
$j_{w*},\, j_{w!}\in Perv_I(\Fl)\subset D_I(\Fl)$.

We have $j_{w_1*}\st j_{w_2*}\cong j_{w_1w_2*}$,
$j_{w_1!}\st j_{w_2!}\cong j_{w_1w_2!}$, $ j_{w_1w_2*}\st j_{w_2^{-1}!}\cong j_{w_1*}$
provided that $\ell(w_1w_2)=\ell(w_1)+\ell(w_2)$.

For $\la\in \La$ the corresponding {\em Wakimoto sheaf} $J_\la\in Perv_I(\Fl)\subset D_I(\Fl)$
is introduced in \cite[3.2]{AB}. It can be characterized by $J_\la\st J_\mu\cong J_{\la+\mu}$ for $\la$,
$\mu\in \La$ and $J_\la\cong j_{\la*}$ for dominant $\lambda$. Notice that $J_\la\cong j_{\la!}$ when
$\la$ is antidominant.

The orbits of $\GO$ on $\Fl$ are indexed by $\La$, we let $\Fl^\la$ denote the orbit corresponding to $\la\in \La$ and let $i_\la:\Fl^\la\to \Fl$ be the embedding.
There exists a unique irreducible Iwahori-Whittaker perverse sheaf on $\Fl^\la$, we let
$\Delta_\la$ (respectively, $\nabla_\la$) be its $!$ (respectively, $*$) extension to $\Fl$.
We have $\Delta_\la=\Delta_0\st j_{w!} $, $\nabla_\la=\Delta_0\st j_{w*}$
if $w\in W\cdot \la$.

We set also $J_\la^{IW}=\Delta_0\st J_\la$. The functor $\F\mapsto \Delta_0\st
\F$ is $t$-exact \cite[Proposition 2]{AB},
so we have $J_\la^{IW}$, $\Delta_\la$, $\nabla_\la\in Perv_{IW}$.

Recall the {\em central sheaves} $Z_\la$, $\la \in \La^+$ of \cite{Gce}.


\subsection{The equivalence of \cite{AB} and its relation to the $t$-structures}  

The main result of \cite{AB} is a construction of an equivalence of triangulated categories
\begin{equation}\label{eqAB}
\Phi_{IW}:D^b(Coh^G(\Nt))\cong D_{IW}(\Fl).
\end{equation}

We recall some properties of the $\Phi_{IW}$ that will be used below.
\begin{equation}\label{PhiO} \Phi_{IW}(\O_\Nt)=\Delta_0, \end{equation}
\begin{equation}\label{PhiOla}\Phi_{IW}(\F\otimes \O_\Nt(\la))\cong \Phi_{IW}(\F)\st J_\la,\end{equation}
\begin{equation}\label{PhiVla} \Phi_{IW}(\F\otimes V_\la)\cong \Phi_{IW}(\F)\st Z_\la,\end{equation}
where $V_\lambda$ denotes the irreducible $G$-module with highest weight $\lambda$.

The following technical statement relating the natural $t$-structures in the two sides of \eqref{eqAB} will play a key role in the proof of the main result.

\begin{Prop}\label{iffIW}
For $\F\in D_{IW}(\Fl)$ the following are equivalent:

i) For all $\la\in \Lambda$, $\F\st J_\la\in D^{\leq 0}(Perv_{IW}(\Fl))$. 

ii) $\F\in \< J_\la^{IW}[d]\ |\ d\geq 0\>$.


iii)  $\Phi_{IW}^{-1}(\F)\in D^{\leq 0}(Coh^G(\Nt))$, where $D^{\leq 0}$ is taken with respect to the
tautological $t$-structure on $D^b(Coh^G(\Nt))$.
\end{Prop}

\proof ii) $\Rightarrow$ iii) follows from  \eqref{PhiO}, \eqref{PhiOla} which imply that
$\Phi_{IW}^{-1}(J_\la^{IW})= \O_{\Nt}(\la)$.

To check that iii) $\Rightarrow$ ii) notice that every equivariant coherent sheaf on a quasi-projective
algebraic variety with a reductive group action is a quotient of a line bundle tensored by a representation of the group. In particular, every object in
$ Coh^G(\Nt)$ is a quotient of a sheaf
of the form $V\otimes \O(\la)$ for some $V\in Rep(G)$. It follows by a standard
argument
that every object $D^{\leq 0}(Coh^G(\Nt))\cap D^b(Coh^G(\Nt))$ is
 a direct summand in an object represented by a finite complex of sheaves of the form $V\otimes \O(\la)$ concentrated in non-positive degrees.
 Thus $D^{\leq 0}(Coh^G(\Nt))\cap D^b(Coh^G(\Nt))=\< V_\nu\otimes \O(\la)[d]\>,$
 $d\geq 0$.
So we see that $\Phi_{IW}^{-1}(\F)\in D^{\leq 0}(Coh^G(\Nt))$ iff $\F\in \< J_\la^{IW}\st Z_\nu \>$.
Since $Z_\nu$ admits a filtration with associated graded being a sum of Wakimoto sheaves \cite[3.6]{AB} and $J_\la^{IW}\st J_\mu\cong
J_{\la+\mu}^{IW}$,
we get that $\< J_\la^{IW}\st Z_\nu \>=\< J_\la^{IW} \>$, which yields the implication.

The implication ii) $\Rightarrow$ i) is clear from $J_\la^{IW}\in Perv_{IW}$, $J_\la^{IW}\st J_\mu\cong
J_{\la+\mu}^{IW}$.

Finally to check that i) $\Rightarrow$ ii) we need an auxiliary statement.

\begin{Lem}\label{lem1}  (see \cite[Lemma 15]{AB})
Given $\F\in D_{IW}(\Fl)$ there exists a finite subset $\fS\subset \Lambda$, such that
for $\mu$, $\la \in \La$ we have $i_\mu^*(\F\st j_{\nu !})=0$ unless $\mu \in \fS+\nu$. \epf
\end{Lem}

Now, to check i) $\Rightarrow$ ii) let $\fS\subset \La$ be constructed as in the Lemma,  and let $\nu
\in \Lambda$ be such that both $\{\nu\}$ and  $\fS+\nu$ are contained in the set of antidominant weights.
Using the standard exact triangles connecting a constructible complex, its $*$ restriction to a closed subset and $!$ extension from the open complement we see that
\begin{equation}\label{FDel}
\F\st j_{\nu!} \in \< \Delta_{\la}[d]\ |\
\la \in \fS+\nu, d\in \Zet \>.
\end{equation}
 Since $\nu$ is antidominant, we have $J_{\nu}=j_{\nu!}$; so condition i)
 says that
 $\F\st j_{\nu!}\in D^{\leq 0}(Perv_{IW}(\Fl))$. Given \eqref{FDel}, this  is equivalent to
$\F\st j_{\nu!} \in \< \Delta_{\la}[d]\ |\
\la \in \fS+\nu, d\geq 0 \>$. Since all $\la$ appearing in the last expression are antidominant, we get
that $\F\st J_\nu \in  \< J_\la^{IW} [d]\ |\ d\geq 0 \>$ which yields ii). \qed

\subsection{The $B_{aff}$ action on $D_{I^0}(\Fl)$}\label{BaffactFl}
The two sided cosets of $I$ in $\LG((t))$ are indexed by $W_{aff}$.
 For each $w\in W_{aff}$ fix a representative $\bw\in IwI$. Let $^wI^0=\bw I^0 \bw^{-1}$
 (it is easy to check that the functors defined in the Lemma below
do not depend on this choice, up to a noncanonical isomorphism).
  Let $Conv_w$ denote the quotient of $I^0\times \Fl$ by the action of $I^0\cap ^wI^0$ given by
$g:(\gamma,x)\mapsto (\gamma g^{-1}, g(x))$. The action map descends to a map $conv_w:Conv_w\to \Fl$.

For $\F\in D_{I^0}(\Fl)$ the complex $\bw_*(\F)$ is equivariant with respect to $I^0\cap ^wI^0$, thus
the complex $\underline{\Ce}\boxtimes \bw_*\F$ on $I\times \Fl$ descends to a canonically defined
complex on $Conv_w$, let us denote it by $\F_w$. The following is standard.

\begin{Lem}
There exists an (obviously unique) action of $B_{aff}$ on $D_{I^0}(\Fl)$, such that
for $w\in W_{aff}$,
$$\tii w: \F \mapsto conv_{w*}(\F_w)[\ell (w)].$$
It satisfies:
$$\tii w^{-1}:\F\mapsto conv_{w^{-1}!}(\F_w)[\ell(w)].$$
\end{Lem}

\subsection{A coherent description of $D_{I^0}(\Fl)$}

A proof of the following result will appear in \cite{tobe}, see also announcement in \cite{ICM}.

Define the {\em Steinberg variety} as $St=\gt\times_\g \Nt$.
Let $Av_{IW}:D_{I^0}(\Fl) \to D_{IW} (\Fl)$ be the {\em averaging} functor,
i.e. the adjoint functor to the embedding of $D_{IW}(\Fl)$ into the category
of constructible complexes on $\Fl$, restricted to $D_{I^0}(\Fl)$.

\begin{Thm}\label{Thmtobe}
There exists an equivalence of triangulated categories
$$\Phi:D^b(Coh^G(St))\cong D_{I^0}(\Fl),$$
satisfying the following properties.

a) $\Phi$ intertwines the $B_{aff}$ action from
section \ref{NCS} with that from \ref{BaffactFl}.

b) $pr_{2*} \circ \Phi^{-1} \cong \Phi_{IW}^{-1} \circ Av_{IW}$.
\end{Thm}

We refer the reader to \cite{ICM} and \cite{FG} for a discussion of some other properties of this equivalence.

\subsection{The "new $t$-structure"}
We are now in the position to derive a partial answer to a question of \cite{FG}.

\begin{Prop}\label{propcomp}
For $\F\in D_{I^0}(\Fl)$ the following are equivalent

i) For all $\la\in \La$ we have $ \F\st J_\la\in D^{\leq 0}(Perv_{I^0}(\Fl))$.

i$'$) There exists $\la_0$ such that $ \F\st J_\la\in D^{\leq 0}(Perv_{I^0}(\Fl))$ for $\la\in \la_0-\La^+$.

ii) $\Phi^{-1}(\F) \in D^{\leq 0}(Coh^G(St))$ where $D^{\leq 0}(Coh^G(St))$ is equipped with the
braid positive normalized $t$-structure for $S=\Nt$.

\end{Prop}

\proof  i) $\Rightarrow$ i$'$) is obvious. To see that i$'$) $\Rightarrow$ i)
notice that any $\la\in \La$  can be written as $\la=\la'+\mu$ where $\mu\in \La^+$ and $\la'\in
\la_0-\La^+$; then $\F\st J_\la = (\F\st J_{\la'})\st j_{\mu*}\in D^{\leq 0}(Perv_{I^0}(\Fl))$
since convolution with $j_{w*}$ is right exact as it amounts to taking direct image under an affine morphism.

It remains to show that i) $\Leftrightarrow$ ii).
Assume that i) holds for a given $\F$. In view of Theorem \ref{brpex}a) and compatibility of the $B_{aff}$ actions,  we need to check that for all $b\in B_{aff}^+$, $pr_*(\Phi^{-1}(b(\F)))\in D^{\leq 0}(Coh(\Nt))$. Since $B_{aff}^+$
acts on $D^b(Perv_{I^0}(\Fl))$ by left exact functors, $b(\F)\st J_\la\in D^{\leq 0}(Perv_{I^0}(\Fl))$
for all $b\in B_{aff}^+$, $\la\in \La$. Thus $Av_{IW}(b(\F))\st J_\la\in D^{\leq 0}(Perv_{IW}(\Fl))$ for all $\la$.
Applying Proposition \ref{iffIW} to  $Av_{IW}(b(\F))\st J_\la$ we see that
$\Phi_{IW}^{-1}(Av_{IW}(b(\F)))=pr_*(b(\Phi^{-1}(\F)))\in D^{\leq 0}(Coh^G(\Nt))$, which gives ii).

The converse statement follows from Proposition \ref{iffIW} and the following
\begin{Prop}
Suppose $\F\in D_{I^0}(\Fl)$ is such that $Av_{IW}(b(\F))\in D^{\leq 0}(Perv_{IW}(\Fl))$ for all $b\in B_{aff}^+$. Then $\F\in D^{\leq 0}(Perv_{I^0}(\Fl))$.
\end{Prop}

\proof
Let $\F$ be as in the statement of the Proposition, and
let $n$ be the smallest integer such that $\F\in D^{\leq n}_{I^0}(\Fl)$.
We need to show that $n\leq 0$.

Let $\F_n$ be the $n$-th perverse cohomology sheaf and $L$ be an irreducible
quotient of $\F_n$. It suffices to show that
 \begin{equation}\label{nenul}
 ^{perv}H^0(Av_{IW}(b(L)))\ne 0
 \end{equation}
 for some $b\in B_{aff}^+$: then using right exactness of the action of $b\in B_{aff}^+$ we see that the composed arrow $b(\F)\to b(\F_n[-n])\to b(L[-n])$
induces a surjection on the $n$-th perverse cohomology sheaf, hence
\eqref{nenul} implies that $n$-th perverse cohomology of $Av_{IW}(b(\F))$
does not vanish, thus $n\leq 0$.

We now check \eqref{nenul}. First we claim that there exists $b\in B_{aff}^+$ such that $Av_{IW}(b(L))\ne 0$.
For this we need the following variation of Lemma \ref{lem1}.

\begin{Lem}
Given $\F\in D_{I^0}(\Fl)$ there exists a finite subset $\fS\subset W_{aff}$, such that
for $w_1$, $w_2 \in W_{aff}$ we have $j_{w_2}^!(\wti_1(\F))=0$ unless $w_2 \in w_1\cdot \fS$.
\end{Lem}

\proof is similar to that of Lemma \ref{lem1} (see \cite[Lemma 15]{AB}). \epf

Now take $b=\tilde{\lambda}$ for a dominant weight $\lambda$. We can assume that
if $w=\mu\cdot w_f\in \fS$ where $\fS$ is as in the Lemma with $\F=L$, $w_f\in W$,
$\mu\in\Lambda$, then $\lambda+\mu$ is strictly dominant. Then each left coset of
$W$ in $W_{aff}$ contains at most one element such that the $!$ restriction of
$\tilde{\lambda} (L)$ to the corresponding $I$ orbit is non-zero (no cancelations in the spectral sequence containing exactly one non-zero entry).
If such an element exists then the corresponding costalk of $Av_{IW}(\tilde{\lambda}(L))$ does not vanish; thus $Av_{IW}(\tilde{\lambda}(L))\ne 0 $ for such $\lambda$.

Choose now $b\in B_{aff}^+$ such that $Av_{IW}(b(L))\ne 0$ and $b$ is an
element of minimal possible length satisfying this property.
Notice that $L$ is $I$-equivariant, and for
an $I$-equivariant complex $L$ we have $\wti(L)\cong j_{w*}\st L$
where $\st$ denotes convolution of $I$-equivariant complexes on $\Fl$.
In particular, when $w=s_\alpha$ is a simple reflection we
have $\tilde{s_\alpha}(L)=j_{s_\alpha*}\st L$,
 $\tilde{s_\alpha}^{-1}(L)=j_{s_\alpha!}\st L$. The perverse sheaves
$ j_{s_\alpha!}$, $j_{s_\alpha*}$ are concentrated on the closure of a one dimensional
$I$-orbit; this closure can be identified with $\Pone$, and we denote it by $\Pone_\alpha$.
 We have an exact sequence of perverse sheaves on $\Pone_\alpha$:
 $$0\to \delta_e \to j_{s_\alpha!}\to j_{s_\alpha*}\to \delta_e\to 0,$$
 where $\delta_e$ denotes the sky-scraper at the zero-dimensional $I$-orbit $\{e\}\subset \Pone_\alpha$.
This exact sequence shows that $\tilde{s_\alpha}(L)$ and $\tilde{s_\alpha}^{-1}(L)$
 are isomorphic in the quotient modulo the thick subcategory generated by $L$.
Let $b=\tilde s_{\alpha_1}\cdots \tilde s_{\alpha_n}$ be a minimal decomposition
of $b$. Then our assumptions on $b$ imply that
$$Av_{IW}(b(L))\cong Av_{IW}(b'(L)),$$
where $b'=\tilde s_{\alpha_1}^{-1}\cdots \tilde s_{\alpha_n}^{-1}$. Since the action
of $b$ is right exact, the action of $b'$ is left exact and $Av_{IW}$ is exact,
we see that $Av_{IW}(b(L))$ is a perverse sheaf; thus the assumption on $b$
implies \eqref{nenul}. \epf

 \begin{Rem} The idea of the proof is partly borrowed from \cite[2.2]{BM}.
\end{Rem}

\begin{Cor}
There exists a $t$-structure $\tau_{\new}$ on $D_{I^0}(\Fl)$ given by: $\F\in D^{\leq 0, \new}$
iff $\F\st J_\la \in D^{\leq 0} (Perv_{I^0}(\Fl))$. The composed equivalence
$$D_{I^0}(\Fl)\underset{\Phi}\isol D^b(Coh^G(St)) \underset{\eqref{equiveqnt}}\isor
D^b(A_\Nt \tmod_{coh}^G)$$
sends $\tau_{\new}$ to the tautological $t$-structure on $D^b(A_\Nt \tmod_{coh}^G)$.
\end{Cor}

\begin{Rem}
Corollary provides a positive partial answer to Question 2.1.3 of \cite{FG}.
In more detail, in {\em loc. cit.} the so called new $t$-structure is defined
on a certain Ind-completion of the bounded derived category of finitely generated $D$-modules
on $\Fl$. Then the question is posed whether this $t$-structure induces one
on the original (not completed) bounded derived category (in a footnote the authors
say they expect a negative answer). The above Corollary gives a positive answer to a weaker question: it shows that the new $t$-structure of \cite{FG} induces one on
the bounded derived category of $I^0$-equivariant finitely generated $D$-modules on $\Fl$.
\end{Rem}
\section{Results of Frenkel-Gaitsgory and a proof of Theorem \ref{main}}\label{CRIT}


\subsection{The functor to modules}
Recall the notion of a nilpotent oper  on a punctured formal disc \cite{FG2} (see also \cite{BDop} for a general introduction to the notion of an oper);
let $\Opn$ be the infinite dimensional scheme parameterizing such opers. By definition
$O\in \Opn$ is a collection of data $O=(\EE,F,\nabla)$ where $\EE$ is a $G$-bundle
on the formal disc $\fD=Spec(\Ce [[t]])$, $\nabla$ is a connection on $\EE$, having a first order pole
at the origin $x_0\in \fD$,  $\nabla: ad(\EE)\to t^{-1} ad(\EE)\Omega^1_\fD$; and $F$ is a $B$-structure on the bundle $ad(\EE)$, these should  satisfy a certain compatibility condition.

The compatibility implies in particular that the residue of the connection
is nilpotent and preserves
the $B$-structure on the fiber at the closed point $x_0$; thus we get a canonical
map $\Opn\to \Nt/G$. We will say that a point $(e,\fb)\in \Nt$ is compatible with a given nilpotent oper if it lies in the corresponding $G$-orbit.

The space of all opers maps isomorphically to the spectrum of center of the category $U_{crit}\gh\tmod$ by \cite{FF}.

The following result is a direct consequence of \cite{FG} compared to Proposition \ref{propcomp}.

We identify $D_{I^0}(\Fl)$ with the category of $I^0$ equivariant critically twisted $D$-modules on $\Fl$,
this is possible by Riemann-Hilbert correspondence, since the critical twisting is integral.
Then we get the derived functor of global sections from $D_{I^0}(\Fl)$ to the derived category of
$U_{crit}\gh$-modules; in fact it lands in the derived category of $I^0$ monodromic modules \cite{FG}. 


Recall that $\BB_e=\{e\}\Ltimes _\g \gt$.

\begin{Thm}\label{FGThm}
%
a) Fix $O\in \Opn$ and  $e\in \N$ so that $Res(O)$ is in the conjugacy class of $e$. Then
there exists an equivalence
$$\Phi_O: DGCoh(\BB_e)\isor D^b(U_{crit}\gh\tmod_O^{I^0}).$$

b) Fix $\tilde e=(e,\fb)\in \Nt$ compatible with $O$.
 For $\F\in  DGCoh(\BB_e)$ 
 we have: $\Phi_O(\F)\in
D^{>0}(U_{crit}\gh\tmod_O^{I^0})$ iff $\Hom _{D^bCoh(St)}(\GG,i_{\tilde e \,*}\F)=0$
for any $\GG\in D^b(Coh^G(St))$ which belongs to $D^{\leq 0}$ with
respect to the braid positive normalized $t$-structure with $S=\Nt$. Here
$i_{\tilde e}$ is the composed map
$\gt \Ltimes _\g \{e\} \isor St \Ltimes_{\Nt} \{\tilde e\} \to St=\gt\times _\g \Nt$.


\end{Thm}

\proof Part (a) is \cite[Corollary 0.6]{FG}.

To check (b) we need to recall the idea of the proof of \cite[Corollary 0.6]{FG}. That  result  is obtained
by combining our equivalence
of Theorem \ref{Thmtobe}
with
the equivalence
of {\em loc. cit.}, Main Theorem 2 between a certain Ind-completion of $D^b(U_{crit}\gh\tmod^{I^0})$ and an (appropriately defined) base change of the
category $D_{I^0}(\Fl)$ with respect to the morphism $\Opn\to \Nt/G$; here
the construction of \cite{AB} is used to endow $D_{I^0}(\Fl)$ with the structure
of a category over $\Nt/G$ (see \cite{FG} and references therein for a definition
of a category over a stack and the notion of base change in this context).
Thus, using the usual notation for base change (and omitting Ind-completion
from notation), \cite[Main Theorem 2]{FG} asserts that:
\begin{equation}\label{MT2}
D^b(U_{crit}\gh\tmod^{I^0})\cong \Opn \times_{\Nt/G} D_{I^0}(\Fl).
\end{equation}

Furthermore, the category $D^b(U_{crit}\gh\tmod_O^{I^0})$
is obtained from $D^b(U_{crit}\gh\tmod^{I^0})$ by base change with respect
to the morphism of point embedding $\{ O\}\to \Opn$. Thus we get
$$D^b(U_{crit}\gh\tmod_O^{I^0})\cong \{O\}\times _{\Opn} (\Opn\times_{\Nt/G}
D_{I^0}(\Fl))\cong \{O\} \times _{\Nt/G}D_{I^0}(\Fl). $$
Substituting the equivalence of Theorem \ref{Thmtobe}:
$D_{I^0}(\Fl)\cong D^b(Coh^G(St))$ we can rewrite the latter category
as
$$\{O\} \times _{\Nt/G}D^b(Coh(St/G))\cong DGCoh(\{O\}\Ltimes _{\Nt/G}
St/G)$$
where the equivalence comes from basic properties of base change for categories. Finally,
$$(St/G)\Ltimes _{\Nt/G}
\{O\} \cong
(\gt/G)\times _{\g/G} (\Nt/G) \Ltimes _{\Nt/G}
\{ O\}
\cong \gt \Ltimes_\g \{O\}
\cong \gt \Ltimes_\g \{e\},$$
which gives the desired equivalence. Notice that the tautological functor
$D^b(U_{crit}\gh\tmod_O^{I^0})\to D^b(U_{crit}\gh\tmod^{I^0})$ corresponds
under the above equivalences to the push-forward functor
$i_{O*}:\{O\}\times_{\Nt/G} D_{I^0}(\Fl)\to \Opn \times_{\Nt/G} D_{I^0}(\Fl)$
which comes from the point embedding $i_O:\{ O\} \to \Opn$ via functoriality
of the base change construction.

We are now ready to deduce statement (b) from
 the
exactness statement in \cite[Main Theorem 2]{FG}. The latter yields the following
description of

 $\Phi_O^{-1}(D^{>0}(U_{crit}\gh\tmod_O^{I^0}))$
 (cf. definition
of the $t$-structure in {\em loc. cit.}, 3.6.1):

for  $\F\in  DGCoh(\{e\}\Ltimes _\g \gt)$ we have
  $\Phi_O(\F)\in
D^{>0}(U_{crit}\gh\tmod_O^{I^0})$ iff the following holds.
For any $\GG\in D_{I^0}(\Fl)$ such that $J_\lambda \st \GG\in D^{\leq 0}(Perv_{I^0}(\Fl))$ for all $\lambda$ we have
\begin{equation}\label{condi}
\Hom_{\Opn \times_{\Nt/G} D_{I^0}(\Fl)} (pr_2^*(\GG), i_{O*}(\F' ))=0,
 \end{equation}
 where $\F'$ is the image of $\F$ under the equivalence 
 $DGCoh(\{e\}\Ltimes _\g\gt)\cong \{O\}\times _{\Nt/G}D_{I^0}(\Fl)$, and  
 $pr_2^*$ denotes the natural pull-back
functor $D_{I^0}(\Fl)\to \Opn \times_{\Nt/G} D_{I^0}(\Fl)$.

Since the composed morphism $pt\to \Nt/G$ is isomorphic to the composition
$\{\tilde e\}\to \Nt\to \Nt/G$, we can use projection formula to simplify the last condition to:
$$\Hom_{DGCoh(\{\tilde e\}\Ltimes_{\Nt} St)}(i_{\tilde e}^*(\Phi^{-1}(\GG)), \F)=0,$$
where $i_{\tilde e}$ is the map $\{\tilde e\}\Ltimes _\Nt St \to St$.
In view of Proposition \ref{propcomp} we get the result. \epf

\subsection{The proof of Theorem \ref{main}}\label{COMP}
In view of Theorem \ref{FGThm} it suffices to check that the image of  $D^{>0}(A_e\tmod^{fg})$ under
equivalence \eqref{ee}
consists of such objects $\F\in DGCoh(\gt_e)$ that $\Hom(\GG,i_{\tilde{e}*}(\F))=0$
 when $\GG\in D^b(Coh^G(St))$ lies in $D^{\leq 0}$ with respect to the normalized braid positive
$t$-structure with $S=\Nt$. This is immediate from Theorem \ref{brpex}. \epf

\begin{Rem}
The equivalences of \cite{FG} are based on existence of a certain infinite
dimensional vector bundle on the space of Miura opers carrying an action of
$\gh$ at the critical level. According to a conjecture of \cite{FG}
the fibers of this bundle are baby
Wakimoto modules for $\gh$ at the critical level.

Likewise, the equivalence of \cite{BMR}, \cite{BM} are based on the existence
of certain vector bundles on the space $\gt$ and its subspaces. In particular,
in \cite{BMR} a certain vector bundle on the formal neighborhood of a Springer fiber
in $\gt$ over a field of positive characteristic $k$ is constructed. It carries an action
of the Lie algebra $\g_k$ and its fibers are baby Verma modules for $\g_k$.

Theorem \ref{main} shows that the two vector bundles are related. In particular, consider
their pull-back to a fixed Springer fiber (defined as either derived, or an ordinary
scheme) where the Springer fiber is embedded in the space of Miura opers by fixing
a point in $\Opn$. Then the bundle of baby Verma modules in positive
characteristic is a sum of indecomposable summands; for almost all  values of $p=char(k)$
each such summand can be lifted
to characteristic zero, and the bundle of $U_{crit}\gh$ modules is the sum of the resulting indecomposable bundles (generally with infinite multiplicities).

It would be interesting to find a more direct explanation of this phenomenon.
\end{Rem}

 \end{document}